\documentclass[12pt]{amsart}
\usepackage{amsthm,amsmath,amssymb,amscd,graphics,enumerate}

{
   \newtheorem{theorem}[subsubsection]{Theorem}
   \newtheorem{proposition}[subsubsection]{Proposition}     
   \newtheorem{lemma}[subsubsection]{Lemma}

   \newtheorem{question}[subsubsection]{Question}
}
{\theoremstyle{definition}

   \newtheorem{example}[subsubsection]{Example}
   \newtheorem{definition}[subsubsection]{Definition}
   \newtheorem{remark}[subsubsection]{Remark}
}
\newcommand{\QQ}{{\mathbb{Q}}}

\newcommand{\ZZ}{{\mathbb{Z}}}

\newcommand{\bmu}{{\boldsymbol{\mu}}}

\newcommand{\cG}{{\mathcal G}}

\newcommand{\cO}{{\mathcal O}}
\newcommand{\cP}{{\mathcal P}}
\newcommand{\cQ}{{\mathcal Q}}

\newcommand{\cX}{{\mathcal X}}
\newcommand{\cY}{{\mathcal Y}}

\newcommand{\Spec}{\operatorname{Spec}}

\newcommand{\End}{{\operatorname{End}}}

\newcommand{\Aut}{{\operatorname{Aut}}}

\newcommand{\smooth}{{\operatorname{sm}}}

\newcommand{\fo}{{\mathfrak{o}}}

\newcommand{\fB}{{\mathfrak{B}}}
\newcommand{\fM}{{\mathfrak{M}}}
\newcommand{\fm}{{\mathfrak{m}}}

\newcommand{\fO}{{\mathfrak{O}}}

\newcommand{\double}{\genfrac..{0pt}1
{\raise -1pt\hbox{$\scriptstyle\longrightarrow$}}{\raise 3pt\hbox
{$\scriptstyle\longrightarrow$}}} 

\def\gen{{\rm gen}}

\def\tototi{\mathbin{\mathop{\otimes}\limits^{\raise-1pt\hbox
{$\scriptscriptstyle {\rm L}$}}}}

\def\indlim{\mathop{\vrule width0pt height7pt depth
4pt\smash{\lim\limits_{\raise 1pt\hbox to 14.5pt
{\rightarrowfill}}}}}
\def\projlim{\mathop{\vrule width0pt height7pt depth
4pt\smash{\lim\limits_{\raise 1pt\hbox to 14.5pt
{\leftarrowfill}}}}}

\begin{document}

\title[Raynaud's group-scheme]{Raynaud's group-scheme and reduction of
coverings}  
\author[D. Abramovich]{Dan Abramovich \\ with an appendix by Jonathan Lubin}
\thanks{Research partially supported by NSF grant DMS-0070970, a
Forscheimer Fellowship and Landau Center Fellowship}  
\address{Department of Mathematics,
Box 1917, 
Brown University, Providence, RI, 02912}
\email{abrmovic@math.brown.edu}
\maketitle

\section{Introduction}
\subsection{Reduction of coverings of degree $p$}
Let $R$ be a discrete valuation ring of mixed characteristics,
$S=\Spec R$, which has generic point $\eta$ with fraction field $K$ and
special point $s$ with residue field $k$ of characteristic
$p>0$. Consider a generically smooth, stable marked curve $Y\to S$
with an action of a finite group $G$
of order divisible by $p$.  Denote $X = Y/G$. We assume that $G$ acts freely on the
complement of the marked points in $Y_\eta$ (it follows that $G$
respects the branches of all nodes of $Y_s$).

In situations where the order of $G$ is prime to the residue
characteristic, the reduced covering $Y_s \to X_s$ is an admissible
$G$-covering, and a nice complete moduli space of admissible
$G$-coverings exists. An extensive literature exists describing that
situation. However, in our case
where the residue characteristic
divides the order of $G$,
interesting phenomena occur (see
e.g. \cite{A-O:Hurwitz}). The situation was
studied by a number of people; we will concern ourselves with results
of Raynaud \cite{Raynaud} and, in a less direct way, 
Henrio \cite{Henrio}. Related work of Saidi \cite{Saidi1, Saidi2,
Saidi3}, Wewers and Bouw
\cite{Wewers1, Wewers2, Wewers3, Bouw, Bouw-Wewers1, Bouw-Wewers2},
Romagny \cite{Romagny} and others 
provides additional inspiration. 

Thus, in our case where $p\,\big|\ |G|$, the covering $Y \to X$ is no longer
generically  
\'etale on each fiber. It is natural to consider some sort of
group-scheme degeneration $\cG\to X$ of $G$, in such a way that $Y$
might be considered something like
an admissible
$\cG$-covering. Raynaud
(\cite{Raynaud}, Proposition 1.2.1)
considered such 
a degeneration locally at the generic points 
of the irreducible components of
$X_s$, in the special case where $|G| = p$; in our first result below we
will work out
its extension to the smooth locus of $X$ (Theorem
\ref{Th:smooth-curve}), and slightly
more general groups, where $p^2\nmid
|G|$ and the $p$-Sylow subgroup of
$G$  is normal. The case where $p^2
\,\big|\ |G|$ remains a question
which I find very interesting (see example \ref{Ex:not-free} for a
negative result in general and remark \ref{Rem:small-ramif} for a
positive result for small ramification).

One still needs to understand the
 structure of $Y \to X$ at the nodes
 of $X_s$ and $Y_s$. 
 Henrio, working $p$-adic analytically, derived algebraic data along $X_s$,
involving numerical invariants and differential forms, which in some sense
classify $Y_s \to X_s$. 
Our second goal in this note is to present a slightly different approach to
such degenerations at a node, modeled on {\em twisted curves},
i.e. curves with algebraic stack
 structures. The point is that, just
 as in \cite{ACV}, the introduction
 of twisted curves
 allows to replace $Y \to X$ by
 something that is more like a
 principal bundle. Unlike the case of
residue characteristics prime to
 $|G|$, the twisted curves will in
 general be Artin stacks rather than
 Deligne--Mumford stacks.

\subsection{Acknowledgements}
Thanks to Angelo Vistoli for help, and to F. Andreatta, A. Corti, A.J. de Jong,
and N. Shepherd-Barron for patient
ears and useful comments. I also heartily thank Jonathan Lubin, who
pointed me in the direction 
of Example \ref{Ex:not-free}, and in particular saved me from
desparate efforts to prove results when $p^2 | \ |G|$.  Matthieu Romagny's computations led to a big improvements in the results obtained - we hope to complete the construction of a moduli space in the near future.

\section{Extensions of groups schemes and their actions in dimension 1
  and 2}

\subsection{Raynaud's group scheme}
 Raynaud (see \cite{Raynaud},
Proposition 1.2.1, see also  Romagny, \cite{Romagny}) considers the following
construction: let $U$ be integral and let $V/U$ 
be a finite flat $G$-invariant morphism of schemes, with $G$ finite. Assume 
that the action of $G$ on the generic fiber of $V/U$ is faithful. We can view
this 
as an action of the constant 
group scheme $G_U$ on $V$, and we consider the schematic image $\cG$ of the
associated 
homomorphism of group schemes 
$$G_U \to \Aut_UV.$$ Since, by definition, $G_U \to U$ is finite, we have that
$\cG 
\to U$ is finite as well.
The scheme $\cG \to U$ can also be recovered as the closure of the image
of the generic fiber of $G_U$, which is, by the faithfulness assumption, 
a subscheme of $\Aut_UV$. By definition $\cG$ acts faithfully on
$V$. 

\begin{definition} We call the scheme $\cG$ the effective model of $G$
acting on $V/U$. 
\end{definition}

Note that a-priori we do not know that $\cG$ is a group-scheme. It is however 
automatically a flat group scheme if $U$ is  the spectrum of a Dedekind
domain. (This follows because, in that case, the image of $\cG \times_U \cG \to
\Aut_UV$ is also flat, and therefore must coincide with $\cG$).

Also note that, if $s$ is a closed point of $U$ whose residue characteristic is
prime to the order of $G$, then the fiber of $\cG$ over $s$ is simply $G$. So
this faithful model is only of interest when the residue characteristic divides
$|G|$.   

I recall and slightly extend a result of Raynaud (see \cite{Raynaud},
Proposition 1.2.1): 

\begin{proposition}\label{Prop:Raynaud}
Let $U$ be the spectrum of a
discrete valuation ring, 
with special point $s$ of residue characteristic $p$ and generic point
$\eta$. Let $V\to U$ be a finite 
and 
flat morphism, and assume that the fiber $V_s$ of
$V$ over $s$ is   
reduced (but not assuming
geometrically reduced).  Assume
given a finite group $G$, with 
normal $p$-Sylow subgroup, {\em
such that $p^2\nmid \ |G|$}, and an
action of $G$ on $V$ such
that $V\to U$  is $G$-invariant, and 
such that the generic fiber $V_\eta \to \{\eta\}$ is
a principal homogeneous space. 
Let $\cG \to U$ be the effective model of $G$
acting on $V/U$. 

 Then $V/U$ is a principal bundle under
the action of  $\cG \to U$.
\end{proposition}

\begin{remark}
An analogous result in a much wider array of cases, in particular relative to a non-finite map, is studied in Romagny's \cite{Romagny-effective}.
\end{remark}

{\bf Proof.} As in Raynaud's
argument, it suffices to show that
the stabilizer of the diagonal of
$V_s \times_UV_s$  inside the group
scheme $V_s \times_U \cG$ is
trivial. Since $G$ acts transitively
on the closed points $t_i$ of $V_s$
sending the stabilizer on $t_i$ to
that over $t_j$, it is enough to
show that one of these stabilizers,
say over $t\in V_s$, is trivial. But
this stabilizer $P$ is a group
scheme over a field with $\deg P | p$, and if nontrivial it
is of degree exactly $p$. In such a
case it must coincide with the
pullback of the 
unique $p$-Sylow group-subscheme of
$\cG$, therefore that $p$-sylow acts
trivially, contradicting the fact
that $\cG$ acts effectively.\qed
\begin{remark} In case the inertia
  group is not normal, Raynaud
  passes to an auxiliary cover,
  which encodes much of the behavior
  of $V \to U$.
\end{remark}
\begin{question}    What can one say on
  the action of $\cG$ on $V$ in 
  case the order of $G$ (and the degree of $V \to U$) is divisible by
   $p^2$, but the inertia
  group is still normal? Specifically, what happens 
  if $|G|=p^2$?
\end{question}
In the latter case, consider a subgroup $P\subset G$ of order $p$. It
can be argued as in Raynaud's proof, that the effective model $\cP\to U$ of
$P$ acts freely on $V$, and thus $V \to V/P$ is a principal
$\cP$-bundle. Similarly, if $\cQ$ is the effective model of $G/P$
acting on $V/P$, then $V/P \to U$ is a principal $\cQ$-bundle. At the
same time, we have an action of the effective model $\cG$ of $G$ on
$V/P$, but it is not necessarily the case that $\cG/\cP \to \cQ$ is an
isomorphism.

As probably the simplest example, consider an action of $\cG_0 =
(\alpha_p)^2=\Spec k[a,b]/(a^p, b^p)$
on $k(t)$.
%
Examples of liftings of a non-free action of the type
$$ t\ \ \mapsto \ \ t\ +\ a\ +\ f(t)\,b$$
for any residue characteristic have been written down by Romagny (personal communication) and Saidi (see \cite{Saidi-example}). The case of
\begin{equation} t\ \ \mapsto \ \ t\ +\ a\ +\ t^p\,b\label{Eq:action}\end{equation}
is particularly appealing, as it involves torsion and endomorphisms of a formal group. I therefore ask

\begin{question}\label{Question}
Can one lift the action (\ref{Eq:action}) to characteristic 0?
\end{question}
A formal positive answer in arbitrary residue characteristics is given by Jonathan Lubin in the appendix. 
Here I consider the case of residue characteristic 2, where  this action can be obtained
as a reduction of an action of $(\ZZ/2\ZZ)^2$ on a smooth curve. I
concentrate on the local picture (making it global is not difficult):\begin{example}\label{Ex:not-free} 
Let $R = \ZZ_2[\sqrt{2}]$. Consider the group-scheme $Y/R$ defined
by
$$ t*t'\ \ =\ \ t\,+\,t'\ +\ \sqrt{2}\,t\, t' $$
This is an additive reduction of the multiplicative group. The
reduction of the subgroup $\mu_2$ is  given as
$$\Spec\, R[a]\, \Big/\, \big(\,a\,(a+\sqrt 2)\,\big),$$ reducing to
$\alpha_2$. It acts on $Y$ by 
translation via the addition law as 
above: $$ t\ \ \mapsto \ \ t\,+\,a\ +\ \sqrt{2}\,a\, t .$$ 

The reduction 
of the action of $\ZZ/2\ZZ$ by inversion is the same group scheme,
again reducing to $\alpha_2$, 
which we   
write as   
$$\Spec\, R[b]\, \Big/\, \big(\,b\,(b+\sqrt 2)\,\big).$$ This
time the action is given by 
$$ t\ \  \mapsto\  \   \left(  1 + {\sqrt{2}}\,b \right) \,t\ \ - \ \
\frac{b\,t^2}{1 + {\sqrt{2}}\,t}. $$ 

Since $2$-torsion is fixed by inversion, these actions
commute. Explicitly, the action of the product is given by 
$$  t  \ \ \mapsto \ \  a\ + \ \left( 1\, +\, {\sqrt{2}}\,b
\right) \,t\ -\  
   \frac{b\,t^2}{1\, +\, {\sqrt{2}}\,t}\ + \  
   {\sqrt{2}}\,a\,\left( \left( 1\, +\, {\sqrt{2}}\,b \right) \,t\ -\ 
      \frac{b\,t^2}{1\, +\, {\sqrt{2}}\,t} \right). $$ 
The reduction modulo $\sqrt{2}$ is given by $$ t\ \ \mapsto \ \ t\ +\
a\ +\ t^2\,b,$$ as required.  
\end{example}  

\begin{remark}\label{Rem:small-ramif} Raynaud's arguments do work when $p^2
\big| \, |G|$ if the $p$-Sylow 
group-scheme of $\cG$ has only \'etale and cyclotomic Jordan--H\"older
factors. This is because, in that case, the are no nonconstant group
subschemes in the reduction. In particular this works whenever the
absolute ramification index over $\ZZ_p$ is $<p$. 
\end{remark}

\subsection{Extension from dimension 1 to dimension 2}

Consider now the case where $\dim U= 2$, and $V/U$ finite flat and $G$-invariant as above. 
Consider the S2-saturation $\cG' \to \cG$ of the effective model $\cG$ of the
$G$ action on $V/U$. We have

\begin{lemma}\label{Lem:flat-is-group}
If $\cG' \to U$ is flat then $\cG'$ is a group-scheme acting on $V$.
\end{lemma}
{\bf Proof.} We claim 
that the 
rational map $\cG'\times_U \cG' \to \cG'$ induced by
multiplication in $Aut_{U}V$ is everywhere defined. Indeed the graph of
this map is finite over $\cG'\times_U \cG'$ and isomorphic to it over
the locus where $\cG'\to \cG$ is an isomorphism, whose complement has codimension
$\geq 2$. The S2 property implies that the graph is isomorphic to
$\cG'\times_U \cG'$, and the map is regular. The same works for the
map defined by the inverse in $Aut_{U}V$. This makes  $\cG'$ a 
group-scheme, and the map $\cG' \to  Aut_{U}V$ into a group-homomorphism.
\qed

This applies, in particular, when $U$ is regular:

\begin{lemma} If $U$ is regular, the S2-saturation $\cG'$ of the effective model
$\cG$ is a finite flat group scheme acting on $V$.
\end{lemma}

{\bf Proof.} Since
$\cG'$ is S2 and 2-dimensional, it is Cohen-Macaulay, and being finite
over the nonsingular scheme $U$, it is finite
{\em and flat} over $U$ (indeed its structure sheaf, being saturated, 
is locally free 
over the nonsingular  2-diemnsional scheme $U$). Pulling back to
$\cG'$ the Cohen-Macaulay flat morphism $\cG' \to U$ we have that
$\cG'\times_U\cG'$ is S2 and flat over $\cG'$ and over $U$. 

The rest follows from Lemma \ref{Lem:flat-is-group}. \qed
 
 When the action on the  generic fiber is free, we have more:

\begin{proposition}\label{Prop:CM-is-principal} Let  $U$ be a Cohen-Macaulay integral scheme with $\dim U =
2$. Let $V 
\to U$ be a $G$ invariant, finite, flat and Cohen--Macaulay morphism,  and assume the
action of $G$ on the 
generic fiber is free. Let $\cG\to U$ be the effective model of the
action. Assume that for every codimension-1 point $\xi$, the  action of
the fiber  $\cG_\xi$ 
on 
$V_\xi$ is free. 

Then
  \begin{enumerate} 
\item $\cG \to U$ is a  flat group-scheme, and 
\item The action of  $\cG$ on $V$ is free. 
\end{enumerate}
\end{proposition}

Note that, by Raynaud's proposition, the assumptions hold when $U =
V/G$ is local of mixed
characteristics $(0,p)$, the fibers 
$V_\xi$ are reduced, the $p$-Sylow
of $G$ is normal and  $p^2 \nmid \ |G|$.

{\bf Proof.} 
Consider the s2-saturation $\cG'$ of $\cG$. Since $V \to U$ is flat and Cohen-Macaulay, the same is true for $V \times_UV \to V$ and for $\cG' \times_UV \to \cG'$. Since $\cG'$ and $V$ are Cohen-Macaulay, we have that $V \times_UV$ and  $\cG' \times_UV$ are Cohen-Macaulay, hence S2, as well. The morphism $\cG' \times_UV \to V \times_UV$ induced by the action $\cG' \to\Aut_UV$ is finite birational and isomorphism in condimension 1. By the S2 property it is an isomorphism. In particular we have that $\cG' \times_UV \to V$ is flat, and since $V \to U$ is faithfully flat we have that $\cG' \to U$ is flat. By Lemma  \ref{Lem:flat-is-group} we have that $\cG'\to U$ is a finite flat group scheme acting on $V$, and the isomorphism $\cG' \times_UV \to V \times_UV$ shows that the action is free, in particular $\cG' \to \cG$ is an isomorphism.
\qed

\section{Curves}
\subsection{The smooth locus}
The main case of interest for us is the following:

Let $R$ be a complete discrete valuation ring of mixed characteristic, with
fraction field $K$ of characteristic 0, residue field $k$ of characteristic
$p>0$, and spectrum $S$. Assume $Y\to S$ is a 
stable curve with smooth generic fiber, $G$ a finite group acting on $Y$ over
$S$, and denote $$X = Y/G.$$ We assume that the closure of fixed points of $G$ in
$Y_K$  are 
disjoint integral 
points of the smooth locus $Y_\smooth$ and hence for every node $y\in Y$, the
stabilizer in $G$ of $y$ keeps the branches of $Y$ at $y$ invariant.   We
denote the complement of the closure of generic branch points 
in $Y_\smooth$ by 
$Y_\gen$, and the image in $X$ by $X_\gen$ - the so called {\em
 general locus}. 

Note that $Y_\smooth \to X_\smooth$ is flat.

The propositions above give:
\begin{theorem}\label{Th:smooth-curve}
Assume  $p^2 \nmid |G|$ and the
$p$-Sylow subgroup of $G$ is normal. 

There exist
\begin{enumerate}
\item a finite flat group scheme $\cG_\smooth\to X_\smooth$, 
\item a homomorphism $G_{X_\smooth} \to \cG_\smooth$ which is an isomorphism on $X_K$,
and 
\item an action of $\cG_\smooth$ on $Y_\smooth$ through which  the action of
$G$ 
factors, 
\end{enumerate}
such that $Y_\gen\to X_\gen$ is a principal $\cG$-bundle.

The formation of $\cG$ commutes with any flat and quasi-finite  base change 
$R \subset R'$. 
\end{theorem}

It would be really interesting to see what happens for other groups $G$. 

\subsection{The structure  of $Y$ and $G$ over nodes of $X$} 

What can be done about the singular points of $X$ and $Y$? It is easy to see
that even in the case of characteristic 0, the cover $Y\to X$ is not a
principal bundle in 
general (and certainly not at the closure of generic branch points). However,
the behavior of $Y\to X$ at the nodes is very 
interesting. Several authors, notably Henrio \cite{Henrio}, have
attached  invariants  to such a node 
which say much about the cover $Y\to X$. It seems desirable to attach a local
geometric structure which underlies such  invariants. My suggested
approach here 
is to follow the method of \cite{AV,modfam,ACV} using twisted curves.

Let us first investigate $Y \to X$.

Consider a node $p\in X$ where is described locally by the equation $xy = \pi^m$, with $\pi$ a
uniformizer in $S$ (splitting the nodes is obtained after a quadratic  base
change and passing to a completion, or strict henselization, on $X$). Similarly take a node $q\in Y$ over $p$ with local equation $st = \pi^n$. Say the local degree of $Y\to  X$ at $p$  is $d$, so without loss of generality we can write $x = s^d\mu$ and $y=t^d\nu$, where $\mu$ and $\nu$ are units on $Y$. Comparing the Cartier divisors of $x,y,s,t$ and $\pi$ on $Y$ we get that $m=dn$, and $\mu\nu=1$. Note that, since $G$ acts transitively on the points of $Y$ lying over $p\in X$, the degree $d$ is independent of the choice of $q$.

\begin{definition} Let
$$Z= \Spec R[u,v]/(uv-\pi^n)$$ on which $\mu_d$ acts via $(u,v) \mapsto (\zeta u , \zeta^{-1} v)$, and consider the $\bmu_d$-invariant $Z \to X$  given by $x = u^d,\, y= v^d$. Define a twisted curve $\cX$, which is a tame separated Artin stack given locally over $p$ by 
$$\cX := [Z/\bmu_d].$$
\end{definition} 

Clearly the coarse moduli space of $\cX$ is $X$. It is not hard to see that $\cX\to X$ only depends on the local ramification degree  $d$ at each node, and not on the choice of local coordinates.

\begin{lemma}
There is a lifting, unique up to a unique isomorphism, of $Y \to X$ to a finite flat Cohen-Macaulay morphism  $Y \to \cX$.
\end{lemma}

{\bf Proof.} 
Consider the $\bmu_d$-cover $P\to Y$ given by $$P =\Spec \cO_Y[w]/(w^d-\mu),$$ where $\mu_d$ acts via $w \mapsto \zeta w$. Define a morphism $P \to Z$  via $u = sw$ and $v = t/w$. This morphism is clearly equivariant, giving a morphism $Y \to [Z/\mu_p]$.  Since $(sw)^p = x = u^p$ and $(t/w)^p = y = v^p$ this lifts the maps to $X$. The uniqueness statement follows from the fact that $\cX$ is separated. To check that $Y\to \cX$ is flat it suffices to show $P \to Z$ flat. This follows from the local criterion for flatness: the fiber over $u=v=0$ is given by $s=t=0, w^d=c$ where $c$ is the constant coefficient of $\mu$ at $s=t=0$. This is a scheme of degree $d$ as required.  Since $Y$ and $\cX$ (or, for that matter, $P$ and $Z$) are Cohen-Macaulay, the morphism is Cohen-Macaulay. \qed

We now have:
\begin{theorem}\label{Th:nodal-curve}
Assume  $p^2 \nmid |G|$ and the
$p$-Sylow subgroup of $G$ is normal. 

There exist
\begin{enumerate}
\item a twisted curve $\cX \to X$, 
\item a finite flat group scheme $\cG\to \cX$, 
\item a homomorphism $G_{\cX} \to \cG$ which is an isomorphism on $\cX_K$,
\item a lifting $Y \to \cX$ of $Y \to X$,
and 
\item an action of $\cG$ on $Y$ through which  the action of
$G$ 
factors, 
\end{enumerate}
such that $Y\to \cX$ is a principal $\cG$-bundle.

The formation of $\cG$ commutes with any flat and quasi-finite  base change 
$R \subset R'$. 
\end{theorem}

{\bf Proof.} There are two issues we need to resolve here: the construction of $\cG$ at the nodes,  and the construction of $\cX$ and $\cG$ at the markings. 

First we need to  extend
$\cG$ over nodes. We  we have that $Y \to \cX$ is  flat and Cohen-Macaulay at the nodes; by Theorem \ref{Th:smooth-curve} $Y_\gen\to X_\gen$  is a principal bundle under $\cG_\gen$.  By  Proposition \ref{Prop:CM-is-principal} the effective model $\cG$ of the action is a finite flat group scheme and $Y$ is a principal bundle.

Next, we deal with the markings: the local picture of $Y_K \to X_K$ at a marking is
$Y=  \Spec R[s]$ and $X = \Spec R[x]$ where $x=s^d$, and the stabilizer in $G$ of $s=0$ on $Y$ is  identified with $\bmu_d$, acting via $s \mapsto \zeta s$. 
Locally around the marking $s=0$ we can define $\cX = [\Spec \cO_X[u]/\bmu_d]$. The discussion above shows that $Y_K \to \cX_K$ is a principal $G$-bundle. Applying Proposition \ref{Prop:CM-is-principal} again we obtain that $\cG$ is a group scheme and $Y$ is a principal bundle. \qed

So, in view of \cite{ACV},  we might
call $\cY \to \tilde\cX$ a
 twisted  $\cG$-bundle.

This suggests an approach to lifting covers from characteristic $p$ to
characteristic 0, by breaking it in  two stages: (1)
lifting group-schemes over $\cX$, and (2)
lifting the covers. Recent work of Wewers \cite{Wewers3} seems to
support such an approach. 

\appendix 
\section{Lifting a non-free action on a formal group}\begin{center}by Jonathan Lubin\footnote{\noindent 
{\sc Department of Mathematics, Brown University Box 1917, Providence, RI 02912}\\ \hspace*{.5cm} {\em Email address:} {\tt lubinj@math.brown.edu}}
\end{center}
\subsection*{The Question} In characteristic $p > 0$, consider the substitution $t\mapsto a + t + bt^p$ , where $a^p = b^p = 0$. This clearly defines a groupscheme of rank $p^2$, isomorphic to $\alpha_p \times \alpha_p$, and an action of the groupscheme on a curve, in this case the affine line.  In Question~\ref{Question} it was asked  whether this group-scheme and this action can be lifted to characteristic zero, over a suitably ramified extension of $\ZZ_p$. 
 
 \subsection*{The Answer} Maybe, depending on what you're willing to accept as an example. I give in this note a lifting of $\alpha_p \times \alpha_p$ to a ring $\fO$, the integers in a well-chosen ramified extension of $\QQ_p$, and an action of this groupscheme on the formal affine line, namely the formal spectrum of $\fO[[t]]$. But if you demand an example of an action on a genuine algebraic curve, I'm afraid that you'll have to look to someone else. 
 
 In general, if $R$ is a ring and $f$ and $g$ are power series in one variable over $R$, then it makes no sense to compose the series, $f\circ g$, unless $g$ has zero constant term. Yet, there are situations where $R$ has a suitable complete topology, when $f\circ g$ can make sense even when $g(0)\neq 0$. Let us detail one fairly general such situation:
 
  If $(\fo, \fm)$ is a complete local ring, then on the category of complete local $\fo$-algebras $(R, M)$ we define a group functor denoted $\fB$ or $\fB_\fo$, such that $\fB(R)$ is the set of power series $$f(t) = \sum_{j\geq 0} c_j t^j \ \ \in \ \  R[[t]]$$ for which $c_0 \in M$ and $c_1 \not\in M$. Our desire is that $\fB_\fo(R)$ should be a group under composition of power series, and indeed the condition on $c_0$ guarantees that composition will be well-defined, while the condition on $c_1$ guarantees that the series will have an inverse in $\fB(R)$. One sees now that if $\kappa$ is the characteristic-p field of definition in Question~\ref{Question}, and if $R$ is the local $\kappa$-algebra $\kappa[a, b]/(a^p, b^p)$, then the series $a + t + bt^p$ is an element of $B_\kappa(R)$. The relation $$(a + t + bt^p)\circ (a' + t + b't^p) = (a + a') + t + (b + b')t^p$$ shows that the groupscheme that's being described is finite and isomorphic to $\alpha_p \times \alpha_p$.

  \subsection*{The Method} We take a formal group $F$ of finite height that has a subgroup of order $p$ as well as a group of automorphisms of order $p$. Now, finite groups of automorphisms of a formal group of finite height are always \'etale, but by taking a slight blowup of $F$, we convert the automorphism subgroup to a local groupscheme, without going so far as to make the above group of torsion points of $F$ \'etale as well. Then, allowing ourselves a slight abuse of language, our desired lifting consists of all substitutions $$t \mapsto a \tilde+ t \tilde+ [b]_{F'} (t),$$ where $a$ is a torsion point of the blow-up of $F$, and $1+b$ is a $p$-th root of $1$. In the displayed formula, $F'$ is the blown-up version of $F$, the tilde over the plus-sign indicates addition with respect to $F'$, and as usual, $[b]_{F '} (t)$ is the endomorphism whose first-degree term is $bt$. I suppose that very confident people may be able to look at the preceding explanation and say, Of Course, No Problem, End of Story. But I'm not so confident, and the rest of this note is devoted to filling in the gaps and making sure, to my own satisfaction at least, that everything is on the up and up. To those confident readers, everything from here on may thus be unnecessary, though the summary 1-7 at the end of this note may be an aid to flagging assurance.

  \subsection{Some Algebra} Let $\zeta = \zeta_p$ be a primitive $p$-th root of 1 in an algebraic extension of $\QQ_p$, and let $\fo = \ZZ_p[\zeta]$. Let also $\pi = \zeta - 1$, a prime element of $\fo$, and let $k$ be the fraction field of $\fo$. In the ring $\fo[T]/(T^p - 1)$, let us call $\Gamma$ the image of $T$, and let us consider $\Delta = \frac{\Gamma-1}{ \pi}$ . Then the minimal polynomial for $\Delta$ is 
  \begin{equation}T^p + \frac{p}{ \pi} T^{p-1} +\frac{p(p - 1)}{ 2\pi^2} T^{p-2} + \cdots+ \frac{p(p - 1)}{ 2\pi^{p-2}} T^2 + \frac{p}{ \pi^{p-1}} T , \tag{$*$}\label{Eq:*}\end{equation}
   in which the coefficient of $T$ is a unit in $\fo$. Let us call $B$ the ring $\fo[\Delta]$; we need to establish a few facts about it. I will use capital Greek letters for elements of $B$, lower case Greek letters for elements of $\fo$. 
   
   \begin{lemma}\label{LubinLemma1} The ring $B$ is isomorphic to $\fo\oplus \fo \oplus\cdots \oplus \fo$, with $p$ factors. In $B$, every element $\Theta$ satisfies the condition that $\Theta^p - \Theta \in \pi B$.
   \end{lemma}
\noindent {\bf Proof.} Recall that $k$ is the fraction field of $\fo$; since the ring $B$ is free over $\fo$ of rank $p$, $B\otimes _\fo k$ is a $k$-vector space of dimension $p$. I will identify a set $\{\Lambda_\xi\}$ of orthogonal idempotents in $B$, indexed by the set of all $p$-th roots of 1, namely 
   $$\Lambda_\xi = \frac{\xi}{ p}\sum^{ p-1}_{i=0} \xi^i\Gamma^{p-1-i} .$$
    Since the summation factors as 
    $$
    \prod_{\rho\neq\xi} (\Gamma - \rho),
    $$ 
     where the product is over $p$-th roots of 1 unequal to $\xi$, and where each factor is the product of $\pi$ with an element of $B$, we see that $\Lambda_\xi$ is indeed in $B$. That the $\Lambda$'s are orthogonal idempotents is easily verified. As a result, the $\Lambda$'s are $\fo$-linearly independent, and therefore they form a $k$-basis of $B \otimes_\fo k$. If we write, for any $z \in B \otimes_\fo k$, it's decomposition $z = \sum_\xi m_\xi(z)\Lambda_\xi$, then each $m_\xi: B \otimes_\fo k \to k$ is a ring morphism, taking elements of the domain that are integral over $\fo$ to elements of $\fo$, since $\fo$ is integrally closed. But $B$ itself is integral over $\fo$, so that the morphisms $m_\xi$ map $B$ to $\fo$. The first conclusion of the Lemma follows. Since each element $\beta \in \fo$ has the property that $\beta^p - \beta \in \pi \fo$, the corresponding fact is true for $B$ as well. It may be of interest to note that this is not true of the subring $\fo[\Gamma]$ of $B$. 
    
    \subsection{Endomorphisms of the fundamental formal group.} We start with the polynomial $f(t) = \pi t + t^p \in \fo[[t]]$, which has associated to it a unique formal group $F(x, y) \in \fo[[x, y]]$ for which $f \in \End_\fo(F)$, as proved in \cite{LT}. The following is hardly surprising:

    \begin{lemma}  For each $\Theta \in B$, there is a unique series $[\Theta]_F (t) \in B[[t]]$ such that $[\Theta]_F ' (0) = \Theta$ and $f \circ [\Theta]_F = [\Theta]_F \circ f$; this series is an element of $\End_B(F)$. 
    \end{lemma}
    
    This may be proved by using either of the halves of Lemma \ref{LubinLemma1}; if one wishes to use the fact that $\Theta^p - \Theta$ is always in $\pi B$, then the proof of the first Lemma in \cite{LT} goes through word for word. 
    
    The endomorphism ring $\End_B(F)$ contains in particular the series $[\Delta]_F$ and $[\Gamma]_F$ ; the $p$-fold iterate of the latter series is the identity series $t$. And since $\Gamma = 1 + \pi\Delta$, our periodic series $[\Gamma]_F (t)$ may also be written as 
    $$
    F\left(\, t, \ \left( [\pi]_F \circ [\Delta]_F\right) (t)\, \right).
    $$
    If $\beta$ is an element of $\operatorname{ker} ([\pi]_F)$, then the series $\tau_\beta(t) = F(t, \beta)$ commutes with both 
    $$
    [\zeta]_F (t)\ \  =\ \  F \left(\,t, \ [\pi]_F (t)\,\right)
    $$
     and 
     $$[\Gamma]_F (t)\ \  =\ \  F \left(\, t, \ \left([\Delta]_F \circ [\pi]_F\right) (t)\,\right).
     $$
      If only $B = \fo[\Delta]$ had not been an \'etale $\fo$-algebra, we could have taken $\ker([\pi]_F)\times \Spec(B)$ as our desired lifting of $\alpha_p \times \alpha_p$. After all, the points of $\ker([\pi]_F)$ are the $\beta$'s mentioned above, and the points of $\Spec(B)$ are essentially the $p$-th roots of unity $\xi$, and the substitution 
      $$
      t\mapsto F\left(\beta, [\xi]_F (t)\right)
      $$
       would be our lifting of the substitution mentioned in the introduction. There is the additional problem that in case $p = 2$, $F$ is of height one and so $\ker[\pi]$ is not a lifting of $\alpha_p$, but the \'etaleness of the other factor is a much bigger obstacle. 
    
    Because of the form of $f(t) = [\pi]_F (t) = \pi t + t^p$, not only $F$ but also all the $B$-endomorphisms $[\Theta]_F$ have the property that the only nonzero terms are in degrees congruent to 1 modulo $p - 1$. Any such series can be written, that is, in the form $\sum_{j\geq 0} H_j$ where each $H_j$ is a form or monomial of degree $1 + j(p - 1)$. For want of a better term, I'll call any series with this last property $(p - 1)$-lacunary.
    
     Now I want to let $\fO$ be any complete local $\fo$-algebra in which $\pi$ is no longer indecomposable, $\pi = \lambda\mu$ where both $\lambda$ and $\mu$ are nonunits. Minimally, one may take $\lambda =\mu=\sqrt\pi$ and $\fO = \fo [\sqrt\pi]$. Or we may let $\fO$ be the ring of integers in any properly ramified algebraic extension $K$ of $k$, and $\lambda$ any element of $K$ with valuation $0 < v(\lambda) < v(\pi) = 1$. Or, generically, we can take $\fO = 
 \fo[[\lambda, \pi/\lambda]]$, a ring that can be described alternatively as $\fo[[\lambda,\mu]]/(\lambda\mu-\pi)$ or as the set of all doubly infinite Laurent series $\sum_{j\in \ZZ} \alpha_j \lambda^j$ in the indeterminate $\lambda$ and with coefficients $\alpha_j \in \fo$ satisfying the additional condition that $j + v(\alpha_j) \geq 0$, where $v$ is the (additive) valuation on $\fo$ and $k$ normalized so that $v(\pi) = 1$.

  If $G$ is a $(p - 1)$-lacunary series in one or more variables, I will call the $\lambda$-blowup of $G$, denoted $G^{(\lambda)}$, the series formed from $G$ in the following way: if 
  $G = \sum_{j\geq 0} H_j$ , each $H_j$ being homogeneous of degree $1 + j(p - 1)$, then $G^{(\lambda)} = \sum_{ j\geq 0} \lambda^j H_j$.

 When we apply the above operation to $F$ and its endomorphisms, here's what happens: $F^{(\lambda)}$ becomes a formal group whose reduction modulo the maximal ideal of $\fO$ is just the additive formal group $x + y$. The maps $\End_\fo(F) \to \End_\fO(F^{(\lambda)})$ and $End_B(F) \to \End_{B\otimes _\fo\fO}(F^{(\lambda)})$ that take $g(t)$ to $g^{(\lambda)}(t)$ are injections. For any $\Theta \in B$, I will write $[\Theta]^{(\lambda)}$ for $[\Theta]_{F^{(\lambda)}} = \left([\Theta]_F\right)^{(\lambda)}$; then since $[\pi]^{(\lambda)}(t) = \pi t + \lambda t^p = \lambda(\mu t + t^p)$, the new formal group $F^{(\lambda)}$ has at least one nontrivial finite subgroup, namely the set of roots of $\mu t + t^p$, under the group law furnished by $F^{(\lambda)}$, and they certainly are the geometric points of $\Spec \left( \fO[t]/(\mu t + t^p)\right)$, but this is not the kernel of $[\pi]^{(\lambda)}$, since the standard construction of kernel in that case leads to something that's not flat. Rather, if we call $g(t) = \mu t + t^p$, then the finite groupscheme we're talking about is the kernel of $g: F^{ (\lambda)} \to F^{(\lambda 2)}$. 
 
 Seeing just how $F^{(\lambda)}$ has, so to speak, an $\alpha_p$-lifting worth of automorphisms is a little trickier and more unusual. Our aim is to show that the automorphism $[\Gamma]^{(\lambda)}(t)$ of $F^{(\lambda)}$ lies in $\fO[\Delta'][[t]]$, where $\Delta' = \lambda\Delta$ has the $\fO$-minimal polynomial 
 \begin{align}
 &T^p + \frac{p\lambda}{\pi} T^{p-1} + \frac{p(p - 1)\lambda^2}{2\pi^2} T^{p-2} +\cdots+ \frac{p(p - 1)\lambda^{p-2}}{2\pi^{p-2}} T^2 + \frac{p\lambda^{p-1}}{ \pi^{p-1}} T \tag{$**$}\label{Eq:**}
 \\
 &= T^p + \frac{p}{\mu}T^{p-1} + \frac{p(p - 1)}{2\mu^2}T^{p-2}+\cdots+ \frac{p(p - 1)}{2\mu^{p-2}} T^2 + \frac{p}{\mu^{ p-1}} T ; \quad \notag  
 \end{align} 
 Note that this polynomial is congruent to $T^p$ modulo the maximal ideal $\fM$ of $\fO$. 
 
 Now recall that $\Gamma = 1 + \Delta\pi$, so that the series $[\Gamma](t)$, which is periodic of period $p$ with respect to substitution of series, whether we're talking about automorphisms of the original $F$ or of the blown-up $F^{(\lambda)}$, can be written 
 $$
 [\Gamma](t)\ \  =\ \  F\left(\,t,\ ([\Delta] \circ [\pi]\,)(t)\right).
 $$
  Since every element of $B$ is an $\fo$-linear combination of $\{1, \Delta,\ldots, \Delta^{p-1}\}$, we may write 
 $$
 [\Delta]_F (t)\  =\  \Delta t +\sum_{j\geq1} C_j t^{j(p-1)+1} \ \ \in\ \  B[[t]],
 $$
  where, as remarked, each coefficient $C_j$ is an $\fo$-linear combination of the powers of $\Delta$, up to $\Delta^{p-1}$. It follows that $[\Delta]^{(\lambda)}$, the corresponding endomorphism of $F$, has the form 
 $$
 [\Delta]_{F^{(\lambda)}} (t) \  =\  \Delta t + \sum_{j\geq1} C_j \lambda^j t^{j(p-1)+1}\ \  \in\ \  \fO[[t]],
 $$ 
  where the $C_j $'s are the same in both displayed formulas. Now, what of $[\Delta]^{(\lambda)} \circ [\pi]^{(\lambda)} = [\Delta]^{(\lambda)}(\pi t + \lambda t^p)$? Making the indicated substitution gives
 \begin{align*}&
  \Delta\ (\pi t + \lambda t^p) \ +\  \sum_{j\geq1} C_j \lambda^j (\pi t + \lambda t^p)^{j(p-1)+1}\\ 
  =\ \ & \Delta'\ (\mu t + t^p) \  +\ \sum C_j \lambda^{jp+1}(\mu t + t^p)^{j(p-1)+1}.  
  \end{align*}
 But now because $C_j \in B = \fo[\Delta]$, we also have $C_j \lambda^{jp+1} \in \lambda \fO[\Delta' ] = \lambda\fO[\lambda\Delta]$, since the $j$'s all are at least 1. This shows that $[\pi\Delta]^{(\lambda)}(t)$ is a power series with coefficients in $\fO[\Delta']$, and indeed, modulo $\fM$, this series is just $\Delta't^p$. Finally, when we add this series and the series $t$ by means of the formal group $F ^{(\lambda)}(x, y) \equiv x + y \ \mod \fM$, the result, namely $[\Gamma] ^{(\lambda)}(t)$, has coefficients in $B' = \fO[\Delta ']$, and is congruent modulo $\fM$ to $t + \Delta 't^p$. One more remark is necessary, the obvious one that if $\mu\alpha + \alpha^p = 0$, then $[\Delta\pi]^{(\lambda)}(\alpha) = 0$ and $[\Gamma] ^{(\lambda)}(\alpha) = \alpha$.

  In summary, this is what we now have: 
  \begin{enumerate}
  \item  The ring $\fo$ is $\ZZ_p[\zeta]$, where $\zeta = \zeta_p$ is a primitive $p$-th root of unity, and we use the prime element $\pi = \zeta - 1$.
  \item  The ring $\fO$ is any suitably ramified extension of $\fo$, the minimal example being $\fO = \fo[\sqrt\pi ]$. This $\fO$ is the ring over which our liftings and action are defined, and we identify in it elements $\lambda, \mu \in \fO$ with $\lambda\mu = \pi$. 
  \item The formal group $F$ over $\fo$ has $\pi t+t^p$ as an endomorphism and thus has $\fo$ as its ring of ``absolute" endomorphisms (over the ring of integers of any algebraic extension field of the fraction field of $\fo$). Allowing for abuse of language, there is a unique $\fo$-subgroupscheme of $ F$ of rank $p$, namely $\ker [\pi]_F = Spec(A)$, where $A = \fo[[t]]/\left([\pi]_F (t)\right)$.
  
  \item The finite $\fo$-algebra $B$ is $\fo[\Delta]$, where the minimal polynomial for $\Delta$ over $\fo$ is given in formula (\ref{Eq:*}). Algebraically, $B$ is $\fo^{ \oplus p}$, and when we call $\Gamma = 1 + \pi\Delta \in B$, we have $\Gamma^p = 1$. The scheme $\Spec(B)$ is a finite \'etale groupscheme of order $p$; the element 
  $\Gamma \in B$ is a ÒgenericÓ $p$-th root of unity, and the operation of the \'etale groupscheme on the formal-affine line is $t \to [\Gamma]_F (t) = F \left(t, ([\pi] \circ [\Delta])(t)\right)$.

  \item We use $\lambda \in \fO$ to form a sort of blowup of $F$, which we call $F^{(\lambda)}$ and which is described on the preceding page. This formal group has the subgroup scheme $\Spec(A ')$, where $A ' = \fO[[t]]/(\mu t + t^p)$, and this groupscheme acts on the formal-affine line by the substitution $t \to F(\lambda)(a, t)$ when $a$ is any root of $\mu t + t^p$. 
  
  \item We define $\Delta ' = \lambda\Delta \in B \otimes_\fo \fO$, and note that $\Delta '$ has the minimal polynomial over $\fO$ given by (\ref{Eq:**}) on the preceding page. Call $B ' = \fO[\Delta ']$. The periodic power series $[\Gamma]^{(\lambda)}(t)$, originally defined to be in $B\otimes_\fo \fO  [[t]]$, actually is in $B '[[t]]$ and as an element of this ring, it becomes $t + \Delta 't^p$ in $B' \otimes_\fO \fO/\fM[[t]]$.
  
  \item Since the series $[\Gamma]^{(\lambda)}(t)$ and the $F^{(\lambda)}(a, t)$ mentioned in (5) commute, we do indeed have a finite groupscheme, namely $\Spec(A ' \otimes_\fO B ')$, acting on the formal-affine line in such a way that over $\fO/\fM$, the action is $t\mapsto a + t + \Delta 't^p$.
  \end{enumerate}

  \end{document}